\documentclass[12pt]{article} %\nofiles [12pt]
\usepackage{amssymb}
\newcommand{\goesto}{\rightarrow}
\newcommand{\cL}{\mbox{$\cal{L}$}}

\newcommand{\cM}{\mbox{$\cal{M}$}}
\newcommand{\cS}{\mbox{$\cal{S}$}}
\newcommand{\PP}{\mathbf P}
\newcommand{\bp}{\mathbf p}

\newcommand{\ml}{\mbox{\scriptsize ML}}
\newcommand{\bic}{\mbox{\scriptsize BIC}}
\newcommand{\mdl}{\mbox{\scriptsize MDL}}
\renewcommand{\hat}{\widehat}
\renewcommand{\setminus}{\smallsetminus}

\newtheorem{theorem}{Theorem}
\newtheorem{lemma}{Lemma}
\newtheorem{remark}{Remark} 
\newtheorem{prop}{Proposition} 
  
\def\squarebox#1{\hbox to #1{\hfill\vbox to #1{\vfill}}}
\newcommand{\qed}{\hspace*{\fill}\vbox{\hrule\hbox{\vrule\squarebox{.667em}\vrule}\hrule}\smallskip}

\def\LL{\ell} \def\u{u}

\begin{document}
   \begin{center}
     {\bf Two new Markov order estimators} 
\\\mbox{}\\ Yuval
Peres\footnote{Department of Statistics,
University of California, Berkeley.
 peres@stat.berkeley.edu.  
Partially supported
by NSF grants
\#DMS-0104073 and
\#DMS-0244479.} and Paul Shields\footnote{Professor Emeritus of
Mathematics, University of Toledo.
paul.shields@utoledo.edu}
         \end{center}
\centerline{\today}
\begin{quotation} \noindent
{\em Abstract.} We present two new methods for estimating the order 
(memory depth) of a
finite alphabet Markov chain from observation of a sample path. 
One method is based on entropy estimation via recurrence times of patterns,
and the other relies on a comparison of empirical conditional probabilities.
The key to both methods is a qualitative change that occurs when a parameter
(a candidate for the order) passes the true order.
We also present extensions to order estimation for
Markov random fields.
 \end{quotation}
 
 \vfill

\bigskip 
 

\newpage 

\section{Introduction} \hspace{\parindent}
\label{intro-sec} Fix a finite set $A$ and 
let $x_m^n$ denote the sequence $x_m,x_{m+1},\ldots, x_n$, where
$x_i\in A$. A stationary, ergodic, $A$-valued process
$X=\{X_n\}$  is  {\em Markov of order $M=0$} if it
is i.i.d., and {\em Markov of order $M > 0$} if $M$ is the
least positive integer such that $\mbox{$\mathbf
P$}(a_{k+1}|a^{k}_{1})= \mbox{$\mathbf
P$}(a_{k+1}|a^{k}_{k-M+1}),$ for all 
  $a_{1}^{k+1}$ such that $ k\geq M$. A {\em consistent
Markov order estimator} is a sequence of functions
$M_n^*\colon A^n\mapsto
\{0,1,\ldots\}, n\geq 1,$ such that for any $M$ and any
 Markov process $X$ of order $M$, 
	$$		\lim_n M^*_n(x_1^n)=
M, \mbox{ a.s. } 
		$$
(Here and throughout, ``a.s.'' always refers to
the  distribution $\PP=\PP_X$ of $X$.)
In this paper we introduce  two new Markov order
estimators. Both use test functions that depend
on the sample size and a candidate $k$ for the order.
The key to our methods is that as $k$ increases,  our test 
functions exhibit a qualitative change of behavior
when $k$ reaches the true order.

Our estimators use the empirical frequencies of overlapping 
blocks, 
\begin{equation} 
		N_n(a_1^k)=
N_n(a_1^k|x_1^n)\stackrel{\mbox{\scriptsize
def}}{=}|\{i\in[0,n-k]\colon x_{i+1}^{i+k}=a_1^k\}|.
		\label{empirical-def}    \end{equation}
  The corresponding empirical probabilities and conditional
probabilities are 
$$
\hat{P}_n(a_1^{k+1})\stackrel{\mbox{\scriptsize def}}{=}{1
\over n-k} N_n(a_1^{k+1}) \; \mbox{ \rm and }
 \hat{P}_n(a_{k+1}|a_1^{k})\stackrel{\mbox{\scriptsize
def}}{=}N_n(a_1^{k+1})/N_{n-1}(a_1^{k}) \,.
$$
 We also define the $k$-step conditional empirical
entropy,
$$
\hat{h}_k(n)\stackrel{\mbox{\scriptsize
def}}{=} -\sum_{a_1^{k+1}}
\hat{P}_n(a_1^{k+1})\log \hat{P}_n(a_{k+1}|a_1^{k})\,.
$$

\bigskip

Our first method, which we call the {\em entropy
estimator method}, compares 
%$k$-step conditional empirical entropy
$\hat{h}_k(n)$
 with the entropy estimator
$[\LL(n)]^{-1}\log n$, 
where  $\LL(n)$ denotes the length of
the longest initial  block in $x_1^n$ that repeats 
in $x_1^n$ (see \cite{Ornstein-Weiss-compression}
and Section 2 below). 

	\begin{theorem} 
$M_n^*(x_1^n)\stackrel{\mbox{\scriptsize def}}{=} \min\{
k\colon\ \hat{h}_k(n)\le [\LL(n)]^{-1}\log n +
2(\log n)^{-1/4}\}$ is a consistent Markov order estimator.
	\label{recurrence-estimator}       \end{theorem}

\bigskip

Our second method, which we call the {\em 
maximal fluctuation method}, is based on the test function
  			\begin{equation}
		\phi_m(x_1^n)\stackrel{\mbox{\scriptsize
def}}{=} \max_{m<k<f(n)} \:
\max_{a_1^k \in A^k} \left[\hat{P}_n(a_{k}|a_1^{k-1}) -
\hat{P}_n(a_{k}|a_{k-m}^{k-1})\right]N_{n-1}(a_1^{k-1}) \, ,
		\label{testfcn-def}    \end{equation}
 where $f(n)=\log\log n$. Define
$ M_n^\#(x_1^n)\stackrel{\mbox{\scriptsize
def}}{=}\min\{m< n-f(n)\colon\
\phi_m(x_1^n)<n^{3/4} \}$; if the set we are minimizing over 
is empty, then we take  $M_n^\#(x_1^n)=n$.

\begin{theorem} 
$ M_n^\#(x_1^n)$ is a consistent Markov order estimator.
	\label{max-flux-thm}       \end{theorem} 
 
A more general form of Theorem~\ref{recurrence-estimator} 
that 
 allows any entropy estimator with a known rate of
convergence is given in Section~\ref{ent-bding-section}.
  An extension of Theorem~\ref{max-flux-thm} to Markov
random fields is given in Section~\ref{rf-sec}.
Connections to other model selection methods
are given in Section~\ref{addenda-sec}.

Careful proofs of Theorems \ref{recurrence-estimator}
and  \ref{max-flux-thm} are given in Sections 
\ref{ent-bding-section} and \ref{max-fluc-section}, 
respectively. For the reader's convenience,
we first present sketches of the proofs.

\bigskip

{\em Sketch of proof of Theorem \ref{recurrence-estimator}:}
 The  $2(\log n)^{-1/4}$ term incorporates the rate of
convergence of $[\LL(n)]^{-1}\log n$ and  that of
$\hat{h}_M(n)$ to their common almost sure limit,
the entropy $H(X)$ of $X$. Thus
 $\hat{h}_M(n)\le [\LL(n)]^{-1}\log n +
2(\log n)^{-1/4}$, eventually a.s., whence  $M_n^*\leq M$ eventually a.s.
On the
other hand, if
$k<M$ then $\hat{h}_k(n)$ converges a.s.\ to the
$k$-step  conditional theoretical entropy
$H_k(X)$, which exceeds $H(X)$, the almost sure
limit of 
$[\LL(n)]^{-1}\log n+2(\log n)^{-1/4}$.
Therefore  $M_n^* \ge M$ eventually a.s.
\bigskip

{\em Sketch of proof of Theorem \ref{max-flux-thm}:}  If
$m<M$ then there exists $a_1^{M+1}$ such that
$\PP(a_{M+1}|a_1^M)> \PP(a_{M+1}|a_1^{M-m+1})$, and hence
$\phi_m(x_1^n)$ grows a.s.\ like
$cn$, for some $c>0$.  Thus $M_n^\#\ge M$  eventually a.s.
On the other hand, classical large
deviations theory shows that for any
$\epsilon >0$, we have $\phi_M(x_1^n)=o(n^{1/2+\epsilon}),$ 
a.s., so $M_n^\#\le M$   eventually a.s.
 
\medskip

 \section{The entropy estimator method.}
\label{ent-bding-section} \hspace{\parindent}
 We first review some elementary facts about entropy, see
\cite{CT} or \cite{Shields-book} for details.    The conditional 
entropy of the next symbol given 
$k$ previous symbols is defined by 
	$$
		H_k=H(X_{k+1}|X_1^{k}) \stackrel{\mbox{\scriptsize
def}}{=}-\sum_{a_1^{k+1}} \mbox{$\mathbf
P$}(a_1^{k+1})\log \mbox{$\mathbf
P$}(a_{k+1}|a_1^{k}).
		$$
 The sequence $\{H_k\}$ is nonincreasing  with limit equal to 
the entropy $H=H(X)$ of the process.   Furthermore, the  process is
Markov of order $M$ if and only if
 \begin{eqnarray*}
    k<M &\Rightarrow& H_k>H(X) \; \mbox{ \rm and }\\
    k\geq M &\Rightarrow& H_k=H(X) \,,
  \end{eqnarray*}
that is, if and only if $H_k$ reaches its limit $H$
{\em exactly} when $k=M$,  see \cite[Thm
I.6.11]{Shields-book}.
  
 The conditional $k$-th order empirical entropy
$\hat{h}_k(n)$ is defined by replacing theoretical
probabilities by the corresponding empirical probabilities.
The ergodic theorem  implies that
for $k$ fixed, $\hat{P}_n(a_1^{k+1})\goesto 
 \mbox{$\PP$}(a_1^{k+1})$ and
$\hat{P}_n(a_{k+1}|a_1^{k})\goesto 
 \mbox{$\PP$}(a_{k+1}|a_1^{k})$, each with
probability 1, and hence
that $
		\hat{h}_k(n)\goesto H_k,$ a.s. Furthermore, in the
Markov case we have the following iterated logarithm result.
	\begin{lemma} If $X$ is Markov of finite order $M$ then
for each $k$ there is a constant  $c_k$ such that 
 		$$
		|H_k-\hat{h}_k(n)|\leq c_k \sqrt{\frac{\log\log n}{n}},
 \mbox{ \rm eventually a.s.\ as } n \to \infty \,.
		$$
\label{LIL}    \end{lemma}
\noindent{\bf Remark.} A slightly weaker inequality
(which would suffice for our application here),
with an extra factor of $\log n$ on the right-hand side can be obtained by 
applying \cite[Theorem 16.3.2]{CT} instead of
(\ref{ent}) below.

  {\em Proof of Lemma \ref{LIL}.} 
Let $\Psi(x)=x\log x-x+1$, so that $\Psi(1)=\Psi'(1)=0$.
 For $x>1/2$ we have $\Psi''(x)=1/x<2$, whence
$|\Psi(x)|<(x-1)^2$ for all $x \ge 1/2$.

Consider two distributions $P$ and $Q$ on the same alphabet $A$, 
and suppose that 
\begin{equation} \label{gam}
\gamma= \max_{a \in A} \Bigl|\frac{P(a)}{Q(a)}-1 \Bigr|  \le 1/2 \, .
\end{equation}
Then the divergence $D(P|Q)=\sum_a P(a)\log\frac{P(a)}{Q(a)}$ satisfies
$$ 
D(P|Q) %\sum_a P(a)\log\frac{P(a)}{Q(a)}
=\sum_a \Bigl[ Q(a) \Psi \Bigl(\frac{P(a)}{Q(a)}\Bigr) +P(a)-Q(a) \Bigr]
=\sum_a \Bigl[ Q(a) \Psi \Bigl(\frac{P(a)}{Q(a)}\Bigr) \Bigr] \leq \gamma^2.
$$ 
Moreover,
$$
\sum_a \Bigl| \Bigl(P(a)-Q(a)\Bigr)\log Q(a)\Bigr| \le
\sum_a \Bigl|\gamma Q(a)\, \log Q(a)\Bigr| = \gamma H(Q) 
\le \gamma \log|A| \,.
$$
Adding the last two inequalities (using positivity of the divergence) gives
\begin{equation} \label{ent}
|H(Q)-H(P)| \le \gamma^2+\gamma \log|A| \,.
\end{equation}
under the assumption (\ref{gam}).

By the law of the iterated logarithm for finite-order Markov
chains, there is a constant $\widetilde{c}_k$ such that
		$$ \biggl|
		\frac{\hat{P}_n(a_1^{k})}{\mbox{
$\PP$}(a_1^{k})} -1 \biggr|\leq
\widetilde{c}_k\sqrt{\frac{\log\log n}{n}}, \; \mbox{ eventually a.s.}
		$$
so an application of (\ref{ent}) to $\hat{P}_n$ and $\PP$ 
proves the lemma.

\qed 

\bigskip

The 
Ornstein-Weiss recurrence theorem,
\cite{Ornstein-Weiss-compression}, states that for any
ergodic finite alphabet process $X$,   the time
until the opening
$n$-block occurs again,
	$$
R_n(x)\stackrel{\mbox{\scriptsize def}}{=}\min\{r\geq n \, \colon \,
x_{r+1}^{r+n}=x_1^n\}, 
		$$
  grows like $e^{nH(X)}$, that is, 
 $(1/n)\log R_n(x)\goesto H(X)$ a.s.    (Earlier, Wyner and Ziv,
\cite{Wyner-Ziv}, established convergence-in-probability
for a related recurrence idea.)  In our setting
$\LL(n)=\max\{k\colon\ R_k\leq n\}$ and the
Ornstein-Weiss recurrence theorem gives
	$$
\lim_{n\goesto\infty}\frac{1}{\LL(n)}\log
R_{\LL(n)}(x)=H(X), \mbox{ a.s.}
		$$
     Let $\cM$ denote the set of ergodic, $A$-valued 
processes $X$ that are finite-order Markov. To obtain a rate
of convergence for $X\in\cM$ we use Kontoyiannis' second-order
result,
\cite[Corollary 1]{Kontoyiannis},  that
for any  $\beta >0$ and $X\in\cM$  
	\begin{equation}
	\log [R_n(x)\mbox{$\PP$}(x_1^n)]=o(n^{\beta}),
\mbox{ a.s.} 
		\label{K-result}    \end{equation}
The statement and proof were for   Wyner-Ziv recurrence but 
can easily be adapted  to Ornstein-Weiss recurrence.
  We use it to prove
\begin{lemma}   $\forall\beta>1/2$ and $X\in\cM$,  $\log
R_n(x)=nH+o(n^{\beta}),$ a.s.
	\label{OW-rate}       \end{lemma} 
{\em Proof.}   Suppose $X$ has order $M$ and
$\beta>1/2$. The Markov
property and the law of the iterated logarithm yield 
 \begin{eqnarray*}
    \log \mbox{$\PP$}(x_1^n)&=&\log
\mbox{$\PP$}(x_1^M)+\sum_{a_1^{M+1}}
N(a_1^{M+1})\log \mbox{$\PP$}(a_{M+1}|a_1^M)\\
&=&(n-M)\sum_{a_1^{M+1}} \mbox{$\mathbf
P$}(a_1^{M+1})\log \mbox{$\mathbf
P$}(a_{M+1}|a_1^M)+o(n^{\beta})\\
&=&-nH+o(n^{\beta}),
 \mbox{ a.s.}
  \end{eqnarray*}
which, combined with (\ref{K-result}), yields
the lemma. \qed
\bigskip
	\begin{lemma}For all $X\in \cM$, 
 
 $\displaystyle 
\;\frac{1}{\LL(n)}\log R_{\LL(n)}\leq
\frac{1}{\LL(n)}\log n\goesto H(X)$, a.s.
\label{logn-form}    \end{lemma}
{\em Proof.} Since $R_{\LL(n)}(x)\leq n\leq R_{\LL(n)+1}(x)$,
the lemma follows from 
		$$
	\frac{1}{\LL(n)}	\log R_{\LL(n)}(x)\leq \frac{1}{\LL(n)}	\log n
\leq  \frac{\LL(n)}{\LL(n)+1}\biggl[	\frac{1}{\LL(n)+1}\log
R_{\LL(n)+1}(x)\biggr],
		$$
  and the fact that both the left-hand and right-hand terms 
go to $H(X)$, a.s. \qed
\bigskip
  
  We also need a lower bound on the growth of $\LL(n)$.
\begin{lemma} For any $X\in \cM$  there is a constant 
$C>0$ such that \newline \hspace*{1em} $\LL(n)\geq
C\log n$, eventually a.s.  \label{bound}
 \end{lemma} 
{\em Proof.} By the Ornstein-Weiss recurrence theorem, 
$$R_k\leq e^{k(H+1)}\leq e^{k(1+\log |A|)}, \; 
\mbox{ eventually a.s. }$$  
 Thus we can take $C=(1+\log |A|)^{-1}$.  \qed
 
 \bigskip

 The lemmas yield 
	\begin{prop} For any $ X\in \cM$,
$$  \frac{1}{\LL(n)}\log n\geq
H(X)-\frac{1}{\sqrt[4]{\log n}}, \mbox{ 
eventually a.s.}$$ 
	\label{recurrence-rate-prop}       \end{prop} 
 {\em Proof.}  The following chain of
inequalities holds \mbox{ eventually a.s. }
		 \begin{eqnarray*}
   \frac{1}{\LL(n)}\log n&\stackrel{\mbox{\scriptsize
(a)}}{\geq}&
		\frac{1}{\LL(n)}\log R_{\LL(n)}
\stackrel{\mbox{\scriptsize (b)}}{\geq} H(X)-
\frac{1}{\LL(n)^{3/8}}\\ &\stackrel{\mbox{\scriptsize
(c)}}{\geq}& H(X)- \frac{1}{[C\log
n]^{3/8}}\stackrel{\mbox{\scriptsize (d)}}{\geq}  H(X)-
\frac{1}{\sqrt[4]{\log n}};		  
  \end{eqnarray*} inequality (a) by Lemma~\ref{logn-form},
 inequality (b) by Lemma~\ref{OW-rate} for $\beta =5/8$
and inequality (c) by Lemma~\ref{bound}, while
inequality (d) is clear.  \qed

\bigskip

We are now ready to prove
Theorem~\ref{recurrence-estimator}, which for ease of
reference we restate here.
  \setcounter{theorem}{0}
  	\begin{theorem} 
$M_n^*(x_1^n)\stackrel{\mbox{\scriptsize def}}{=} \min\{
k\colon\ \hat{h}_k(n)\le [\LL(n)]^{-1}\log n +
2(\log n)^{-1/4}\}$ is a consistent Markov order estimator.
	\label{recurrence-estimator2}       \end{theorem} 
 {\em Proof.} Suppose $X\in \cM$ has order $M$ and entropy
$H=H(X)$.  We first show that underestimation  does not
occur, 
 eventually a.s. For
$m<M$,
 the simple facts
    \begin{enumerate}
         \item[\hspace*{1em}(a)]  $
\hat{h}_m(n)\stackrel{\mbox{\scriptsize a.s.}}{\goesto}  
H_m$ as $ n \to \infty$ and $H_m >H,$
       \item[\hspace*{1em}(b)] $[\LL(n)]^{-1}\log n
\stackrel{\mbox{\scriptsize a.s.}}{\goesto} H$ and
$2(\log n)^{-1/4} \goesto 0,$
   \end{enumerate} 
 immediately imply that 
$\hat{h}_m(n)>  [\LL(n)]^{-1}\log n + 2(\log n)^{-1/4},$ 
  eventually a.s.

\bigskip

  The following chain of
inequalities holds  eventually a.s.
 \begin{eqnarray*}
    \hat{h}_M(n) &\stackrel{\mbox{\scriptsize (a)}}{\leq}&
H+ c\sqrt{\frac{\log\log
n}{n}}\stackrel{\mbox{\scriptsize (b)}}{\leq} 
\frac{1}{\LL(n)}\log n+\frac{1}{\sqrt[4]{\log n}} +
c\sqrt{\frac{\log\log n}{n}}\\ 
&\stackrel{\mbox{\scriptsize (c)}}{\leq}&
\frac{1}{\LL(n)}\log n+\frac{2}{\sqrt[4]{\log n}},
  \end{eqnarray*}
  inequality (a) by Lemma~\ref{LIL} and the fact that
$H_M=H$, and inequality (b) by
Proposition~\ref{recurrence-rate-prop}, while inequality (c)
is obvious. We conclude that
$\cM_n^*(x_1^n)\leq M$,  eventually
a.s.  
\qed

	\bigskip

As the above proof suggests, the entropy
estimator $[\LL(n)]^{-1}\log n$ can be replaced by any
consistent entropy estimator $\widehat{H}(x_1^n)$ that has an
{\em $o(1)$ underestimation bound}, i.e., a
function
$\u(n)\goesto 0$ such that for all $X\in \cM$
		\begin{equation}
		\widehat{H}(x_1^n)\geq H(X)-\u(n), \mbox{ eventually a.s., } 
		\label{good-estimator}    \end{equation}
provided we replace  
$2/\sqrt[4]{\log n}$ by  $|\u(n)|+ (1/n)\log n$.  
\setcounter{theorem}{0}
	\begin{theorem}[General form] \mbox{}
 
 Let $\widehat{H}(x_1^n)$ be a
consistent entropy estimator with $o(1)$ underestimation
bound
$\u(n)$.  Then
		$\displaystyle 
			M_n^*(x_1^n) \stackrel{\mbox{\scriptsize def}}{=}\min\{k\colon\
\hat{h}_k(n)<\widehat{H}(x_1^n)+|\u(n)|+(1/n)\log n\}
$
is a consistent Markov order estimator.
\label{entropy-method}       \end{theorem}

   We used the recurrence-based entropy estimator as it is one of
the simplest to describe and compute, it easily
updates as
$n$ increases, and its second order properties are easy
to determine.  Its underestimation bound
$1/\sqrt[4]{\log n}$ goes to 0 very slowly, however,
which suggests that  its associated order estimator
$M_n^*(x_1^n)$ converges slowly to
$M$.  Furthermore, though the recurrence idea does generalize to
higher dimensions, see
\cite{Ornstein-Weiss-2dim}, a useful rate theory for it has not 
been  established.  In
Section~\ref{other-entropy-sec},  we present another
entropy estimator
that has a more rapidly convergent underestimation bound
and is extendable to higher dimensions.

\section{The maximal fluctuation  method.} 
\label{max-fluc-section} \hspace{\parindent}   
We now prove the second theorem stated in the
introduction, namely, 
  \begin{theorem} 
$ M_n^\#(x_1^n)\stackrel{\mbox{\scriptsize
def}}{=}\min\{m< n-f(n)\colon\
\phi_m(x_1^n)<n^{3/4}
\} $
is a consistent Markov order estimator.
(Recall that we defined $ M_n^\#(x_1^n)=n$ if this set is empty).
\label{max-flux-thm2} \end{theorem} 

{\em Proof.}  
Let
$$ \label{delt}
		\delta_m(a_1^k|x_1^n)
\stackrel{\mbox{\scriptsize def}}{=} 
N_n(a_1^{k})- N_{n-1}(a_1^{k-1}) \hat{P}_n(a_k |a^{k-1}_{k-m}) 
$$
and note that
\begin{equation} \label{phi}
		\phi_m(x_1^n)
=\max_{m<k<f(n)}
\max_{a_1^k}\delta_m(a_1^k|x_1^n) \,,
\end{equation}	
where $f(n)=\log\log n$.
We first show that
eventually a.s.\ underestimation does not occur. Suppose $X\in \cM$
has order $M$ and $m<M$.  Choose $a_1^{M+1}$ such that 
$$
\mbox{$\PP$}(a_{M+1}|a_1^{M})> \mbox{$\mathbf
  P$}(a_{M+1}|a_{M-m+1}^{M}). 
$$
By the ergodic theorem there exists
$\epsilon>0$ such that, eventually a.s., 
$$
N_{n-1}(a_1^{M})>\epsilon n
\mbox{ \rm and } \hat{P}_n(a_{M+1}|a_1^{M})- \hat{P}_n(a_{M+1}|a_{M-m+1}^{M})
\geq \epsilon.
$$
This implies that $\phi_m(x_1^n)\geq \epsilon^2 n$ and hence that
$M_n^\#(x_1^n)\geq M$, eventually a.s.

 It takes somewhat more effort to show that, eventually
a.s., $\phi_{M}(x_1^n)\leq n^{3/4}$. We first 
	note that for fixed
$k\geq M$,
	\begin{equation}
		Z_k(n)\stackrel{\mbox{\scriptsize
def}}{=}N_n(a_1^k|x_1^n)-N_{n-1}(a_1^{k-1}|x_1^{n-1})
\mbox{$\PP$}(a_k|a_{k-M}^{k-1}), \; 
 n\geq k,
		\label{diff-sequence}    \end{equation} is a
martingale with bounded differences.  Indeed, with
$\chi(B)$ denoting the indicator of   $B$, we can write
  $ Z_k(n)=\sum_{j=k}^n\Delta_k(j),
$ where $$
\Delta_k(j)=\chi(X_{j-k+1}^j=a_1^k)-
\chi(X_{j-k+1}^{j-1}=a_1^{k-1})
\mbox{$\PP$}(a_k|a_{k-M}^{k-1}),
$$ 
and direct calculation shows that 
$E(\Delta_k(j)|X_1^{j-1})=0$ and
$ \|\Delta_{k}(j)\|_{\infty} \le 1$  for  $j> k.$	
 {}From the Hoeffding-Azuma large deviations bound
for martingales with bounded differences,
\cite{Hoeffding,Azuma},  the probability that
$|Z_n| \ge n^{3/4}$ is at most
$2\exp(-n^{1/2}/2)$.

A similar argument also shows that for
		$$
	Z_k^*(n)\stackrel{\mbox{\scriptsize def}}{=}
N_n(a_{k-M}^k|x_1^n)-N_{n-1}(a_{k-M}^{k-1}|x_1^{n-1})
\mbox{$\PP$}(a_k|a_{k-M}^{k-1}), \; 
 n\geq k,
		$$
the probability that $|Z_k^*(n)| \ge n^{3/4}$ is at most
$2\exp(-n^{1/2}/2)$. 

Next we note that 
		$$
		Z_k(n)-\delta_M(a_1^k|x_1^n)= \frac{N_{n-1}(a_1^{k-1})}
{N_{n-1}(a^{k-1}_{k-M})}\;
 Z_k^*(n) \,,
		$$
which has absolute value at most $|Z^*_k(n)|$. 
 Thus, the probability that
$\delta_M(a_1^k|x_1^n) \ge 2 n^{3/4}$ is less than
$4\exp(-n^{1/2}/2)$. 
 Since there are at most
$|A|^{f(n)+1}=n^{o(1)}$ possible sequences
$a_1^k$, it follows from (\ref{phi}) and 
an application of Borel-Cantelli that
eventually a.s.,
$\phi_{M}(x_1^n) \leq n^{3/4}$. This completes the proof of
Theorem~\ref{max-flux-thm}. \qed

\begin{remark} 
After one of us lectured on these results~\cite{Shields-talk},
B. Weiss noted that in recent joint work he did with G. Morvai,
they independently developped the estimator
$M_n^\#$ discussed in Theorem \ref{max-flux-thm}. 
\label{weiss-remark}       \end{remark}

\subsection{Markov Random Fields}\hspace{\parindent} 
\label{rf-sec}  The method of maximum fluctuations extends
in  modified form to Markov random fields, where order is
usually called {\em range.}  We confine our discussion to
the two dimensional (2-d) case; the extension to higher
dimensions is straightforward.  

We  use the following notation.
   \begin{enumerate}
 \item $\cS_t\stackrel{\mbox{\scriptsize
def}}{=}\{(i,j)\colon\ -t\leq i\leq t, \;  -t\leq j\leq
t\}=$ the square of width $2t+1$, centered at the origin. 
(Note that $\cS_{t+s}\setminus \cS_t$ is a
square ``annulus'' of thickness $s$.)
 \item $\cS_t(\bar{u})\stackrel{\mbox{\scriptsize
def}}{=}$ the square of width $2t+1$ with center at $\bar{u}\in
\mbox{$\mathbf
Z$}^2$. 
 \item   $\Lambda_n\stackrel{\mbox{\scriptsize def}}{=}$
the square of width $n$ with lower left corner at $(1,1)$.
 
 \item A configuration $a(\Lambda)$  is a function
$a\colon\Lambda\mapsto A$; if no confusion 
results its  restriction to $\Lambda'\subseteq
\Lambda$ will be denoted by
$a(\Lambda')$.
   \end{enumerate}

  A random field is
a collection $X=\{X(\bar{n})\colon\ \bar{n}\in
\mbox{$\mathbb Z$}^2\}$ 
of random variables with values in $A$.
Unless stated otherwise, random fields are
assumed to be stationary and ergodic.  We 
use the conditional probability notation 
$$
\mbox{$\mathbf
P$}(a(\Lambda)|b(\Lambda'))
\stackrel{\mbox{\scriptsize def}}{=}
 \frac{\mbox{Prob}(X(\Lambda)=a(\Lambda),  X(\Lambda'))
 =b(\Lambda'))} 
 {\mbox{Prob}(X(\Lambda')=b(\Lambda'))}.
 $$
A random field  is said to be {\em Markov with range $R=0$}
if it is i.i.d, and {\em Markov with range $R\geq 1$} if  
$R$ is the least positive integer $r$ such that for all $\ell\geq0$ and $
t\geq 0$
	$$
		\mbox{$\mathbf
P$}(a(\cS_{\ell})|b(\cS_{\ell+r+t}\setminus
\cS_{\ell})) =
\mbox{$\mathbf
P$}(a(\cS_{\ell})|b(\cS_{\ell+r}\setminus
\cS_{\ell})), 
		$$
for all configurations $a(\cS_{\ell})$ and
$b(\cS_{\ell+r+t}\setminus \cS_{\ell})$. That is, $R$ is the
least
$r$ such that the random variables
$X(\cS_{\ell})$ on the inner square and 
$X(\cS_{{\ell}+r+t}\setminus \cS_{{\ell}+r})$ on the outer
annulus are conditionally independent, given the values
$X(\cS_{\ell+r}\setminus
\cS_{\ell})$ on the inner annulus.  The range of a
finite-range random field
$X$ is denoted by $R=R(X)$.

\bigskip

Our 2-d maximum fluctuation method tests whether configurations on a square 
are conditionally independent of those  outside a square  that is
expanded by 
$r$ in each axis direction, given the
configuration in the annulus between the two
squares. Not only do we need to test over a
(slowly) growing interval of  possible orders
$r$, but now we also need to examine a (slowly) 
growing interval of 
 sizes $\ell$ for the inner square, as order can
depend on square size, though it  eventually
becomes constant as square size increases.  Counting
overlapping blocks as in (\ref{empirical-def}) will not be
used because the higher dimensional analogue of
(\ref{diff-sequence}) need not be a martingale.
We focus instead  on counting  nonoverlapping
blocks, to which classical large deviations is
applicable, but now we must also consider translates.  

Given  $n>8$, let $\ell,r,$ and $ t$ be integers in the closed
interval
$[0,\log \log n]$  and put 
$$
k\stackrel{\mbox{\scriptsize def}}{=}
\ell+ r+ t, \mbox{ and }  T\stackrel{\mbox{\scriptsize
def}}{=}\biggl\lceil\frac{n}{2k+1} \biggr\rceil -1.$$ 
We assume the integer  $n$ is
large enough to guarantee that  $T>0$ for all 
 $k\leq 3 \log\log n$.  Let
		$$
		\Pi_k=\{\cS_k(\bar{u}_1), \cS_k(\bar{u}_2), \ldots ,
\cS_k(\bar{u}_{T^2})\}
		$$
be the partition of the square $\Lambda_{(2k+1)T}$ into
squares of width
$2k+1$.  For each $\bar{v}\in \Lambda_{2k+1}$, let
$\Pi_k(\bar{v})=\{\cS_k(\bar{v}+\bar{u}_j), 1\leq j
\leq T^2\}$ be the translated partition of the square
$\bar{v}+\Lambda_{(2k+1)T}\subseteq\Lambda_n$. \newline
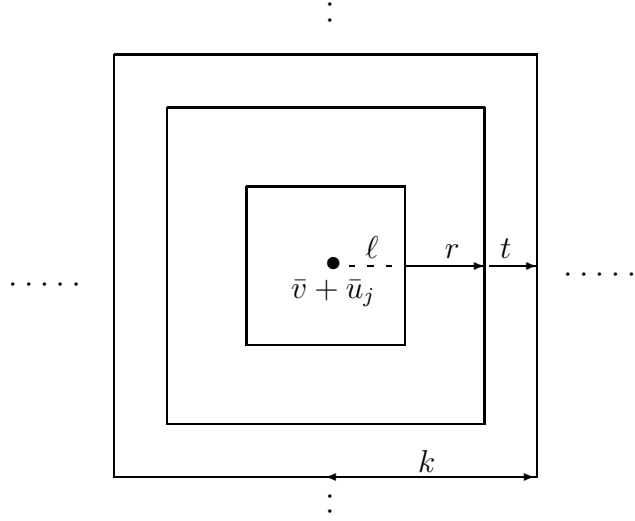
{ \begin{picture}(80,195)(20,60)
%Large box
\put(140,230){\line(1,0){160}}  %box top
\put(140,230){\line(0,-1){160}}  %box left
\put(140,70){\line(1,0){160}}  %box bottom
\put(300,230){\line(0,-1){160}}  %box right
%middle box
\put(160,210){\line(1,0){120}}  % box top
\put(160,210){\line(0,-1){120}}  % box left
\put(160,90){\line(1,0){120}}  % box bottom
\put(280,210){\line(0,-1){120}}  % box right
%inner box
\put(190,180){\line(1,0){60}}  % box top
\put(190,180){\line(0,-1){60}}  % box left
\put(190,120){\line(1,0){60}}  % box bottom
\put(250,180){\line(0,-1){60}}  % box right

\put(220,148){$\bullet$}
\put(207,138){$\bar{v}+\bar{u}_j$}
\put(235,153){$\ell$}
\put(265,153){$r$}
\put(252,150){\vector(1,0){28}}
\put(282,150){\vector(1,0){18}}
\put(286,153){$t$}
\multiput(222,150)(7,0){5}{\line(1,0){2}}

\multiput(310,144)(6,0){5}{$\cdot$}
\multiput(100,140)(6,0){5}{$\cdot$}
\multiput(220,240)(0,6){2}{$\cdot$}
\multiput(220,60)(0,-6){2}{$\cdot$}
%\put(220,69){$\tiny \bullet$} 
\put(227,70){\vector(1,0){72}}
\put(295,70){\vector(-1,0){75}}
\put(255,72){$ k$}
\end{picture}} 

\hspace*{-5em}

Given  a configuration  $x(\Lambda_n)$ and a configuration
$a(\Lambda)$ on  a centrally symmetric subset
$\Lambda\subset
\cS_k$, and given a vector
$\bar{v}\in
\Lambda_{2k+1}$, put 
	 $$
  N_{\bar{v}}(a(\Lambda))=
N_{\bar{v}}(a(\Lambda)|x(\Lambda_n))
\stackrel{\mbox{\scriptsize
def}}{=}\#\{j\colon\
x(\bar{v}+\bar{u}_j+\bar{w})=a(\bar{w}), \forall
\bar{w}\in \Lambda\},
  $$
  that is, the number of times the configuration $a(\Lambda)$
appears in $x(\cdot)$, centered at a member of the translated partition
$\Pi_k(\bar{v})$.  Our 2-d test function is
	\begin{equation}
		 	\delta_{\ell,r,t,\bar{v}}
(a(\cS_{k})|x(\Lambda_n))\stackrel{\mbox{\scriptsize def}}{=}
N_{\bar{v}}(a(\cS_{k}))- N_{\bar{v}}(a(\cS_{k}\setminus
\cS_{\ell}))
\frac{N_{\bar{v}}(a(\cS_{\ell+r}))}{N_{\bar{v}}(a(\cS_{\ell+r}
\setminus\cS_{\ell}))}.
		\label{2d-test-fcn}    \end{equation} 
This is maximized over configurations $a(\cS_{k})$ and
translates
$\bar{v}$ to produce
		$$
		\delta_{\ell,r,t}(x(\Lambda_n)) \stackrel{\mbox{\scriptsize
def}}{=} \max_{\bar{v}\in
\Lambda_{2k+1}}\; \max_{a(\cS_{k})}\;\delta_{\ell,r,t,\bar{v}}
(a(\cS_{k})|x(\Lambda_n))\,.
		$$
  For $\ell=\lfloor \log \log n\rfloor$ define
		$$
		\phi_{r}(x(\Lambda_n)) \stackrel{\mbox{\scriptsize
def}}{=}\max_{0<t<
\log\log n}\;
\delta_{\ell,r,t}(x(\Lambda_n)).	
$$

  Our 2-d order estimator is
  $$ R^*_n(x(\Lambda_n))\stackrel{\mbox{\scriptsize
def}}{=}\min\{r< n-3\log\log n\colon\
\phi_r\leq n^{3/2}\}, 	$$
where, if there is no such
$r<n-3\log\log n$, we set  
$R^*_n(x(\Lambda_n))=n$.

	\begin{theorem} \mbox{}
Let $X$ be a stationary, ergodic, finite range random field on ${\mathbb  Z}^d$.
Then $R^*_n(x(\Lambda_n))=R(X)$, eventually a.s.
	\label{consistency-2d}       \end{theorem} 
{\em Proof.} If $r<R=R(X)$, an argument similar to the
1-dimensional case shows that $\phi_r(x(\Lambda_n))\geq C n^2$,
eventually a.s., for a some  $C>0$.  Thus,  
underestimation eventually a.s.\ does not occur.

To complete the proof it is enough to show that
$\phi_R<n^{3/2}$, eventually a.s. Towards this end,
we fix $\ell>0$ and $ t>0$,  put $r=R$ and
$k=\ell+R+t$, fix
$a(\cS_{k})$ and
$\bar{v}\in
\Lambda_{2k+1}$, and put
$N=N_{\bar{v}}$.  Our 2-d test function
(\ref{2d-test-fcn}) can then be expressed as the sum
 $$
N(a(\cS_{k}))- N(a(\cS_{k}\setminus
\cS_{\ell}))
\frac{N(a(\cS_{\ell+R}))}{N(a(\cS_{\ell+R}
\setminus\cS_{\ell}))}=\Delta_1+\Delta_2,$$
where
 {\small \begin{eqnarray}
    \Delta_1&=&N(a(\cS_{k}))- N(a(\cS_{k}\setminus
\cS_{\ell}))
\mbox{$\mathbf
P$}(a(\cS_{\ell})|a(\cS_{\ell+R}\setminus\cS_{\ell})),\nonumber\\
\mbox{ and}\nonumber\\
\Delta_2&=&
\frac{N(a(\cS_k\setminus\cS_{\ell}))}{N(a(\cS_{\ell+R}
\setminus\cS_{\ell}))}\biggl[
N(a(\cS_{\ell+R}\setminus
\cS_{\ell}))\mbox{$\mathbf
P$}(a(\cS_{\ell})|a(\cS_{\ell+R}\setminus\cS_{\ell}))-
N(a(\cS_{\ell+R}))\biggr]\nonumber
\\\mbox{}\nonumber\\&\leq & \biggl|N(a(\cS_{\ell+R}\setminus
\cS_{\ell}))\mbox{$\mathbf
P$}(a(\cS_{\ell})|a(\cS_{\ell+R}\setminus\cS_{\ell}))-
N(a(\cS_{\ell+R}))\biggr|. \label{delta2-bound}
  \end{eqnarray}}

Denote
$\widetilde{\mbox{$\bp$}}=\mbox{$\mathbf
P$}(a(\cS_{\ell})|a(\cS_{\ell+R}\setminus\cS_{\ell}))$ and
$\bar{w}_j=\bar{u}_j+\bar{v}$.
Then we can write $\Delta_1=\sum_{j=1}^{T^2}\Delta_{1,j}$,
where with $\chi(\cdot)$ denoting the indicator function,
$$
		\Delta_{1,j}\stackrel{\mbox{\scriptsize def}}{=} 
\chi\biggl(X(\cS_k(\bar{w}_j))=a(\cS_k)\biggr)-
\chi\biggl(X([\cS_k\setminus\cS_{\ell}](\bar{w}_j))
=a(\cS_k\setminus\cS_{\ell})\biggr)
\widetilde{\mbox{$\bp$}}.
		$$

 Therefore, conditioned on the  values
$a(\cS_k\setminus\cS_{\ell})$ in the square annulus,
$\Delta_{1}$ is a sum of
$N(a(\cS_{k}\setminus
\cS_{\ell})) \le T^2$ binary i.i.d.\ mean 0 random variables.  The
classical Hoeffding large deviations bound,
\cite{Hoeffding},  implies that the probability that
$|\Delta_1|>\frac{1}{2}n^{3/2}$ is at most
$2 \exp(-n/4)$.  The inequality (\ref{delta2-bound})
implies that the same  result holds for $|\Delta_2|$.  
Since there are only subexponentially
many $a(\cS_{k}\setminus \cS_{\ell})$ and $\bar{v}\in
\Lambda_{2k+1}$ to consider, the Borel-Cantelli lemma
implies that $\phi_R\leq n^{3/2}$, eventually a.s.
This completes the proof of Theorem~\ref{consistency-2d}. \qed

 \begin{remark} 
To simplify the discussion we focused on squares rather than
diamonds which are more natural in  Ising models.   Our
concepts and results can easily be converted to the latter
setting.  
	\label{Markov-def-remark}       \end{remark}  
 
 \begin{remark}  
  Csisz\'ar and Talata, \cite{Csiszar-Talata}, have
recently shown the existence of a consistent range estimator
for a restricted class of Markov random fields, namely,
those for which, conditioned on any boundary,
probabilities in a square are positive, a condition
that allows them to focus only on squares of size 1,
rather than squares of growing size as we did.  They
assume no  bound on the range and use a variant
of the BIC in which maximum likelihood is replaced by
maximum pseudolikelihood.  

\label{CT}       \end{remark}

\section{Extensions and related
work.}%\hspace{\parindent}
\label{addenda-sec}
\subsection{Other entropy
estimators.}\hspace{\parindent} \label{other-entropy-sec}
There are many known consistent entropy estimators,
for most of which $o(1)$ underestimation bounds for the 
Markov case have not been established.  In addition to the
recurrence estimator such an underestimation bound can be
shown to hold for  $\hat{h}_{f(n)}(n)$, where for example 
$f(n)=\frac{\log \log n}{\log |A|}$. 
\begin{prop} There a positive constant $C$ such that for any
$X\in
\cM$, 
    \begin{enumerate}
    \item[\hspace*{1em}(a)] $\hat{h}_{f(n)}(n)\goesto
H(X)$, a.s.
  \item[\hspace*{1em}(b)] $\hat{h}_{f(n)}(n)\geq
H(X)-C\frac{\log^2 n}{n}$, eventually a.s.
   \end{enumerate}
   \label{coding-rate-thm}       \end{prop} 
{\em Proof.}\  By the
Ornstein-Weiss entropy estimation theorem,
\cite{Ornstein-Weiss-sampling}, the per-symbol  empirical block
entropy 
$\frac{1}{f(n)}H(\hat{P}_{f(n)}(\cdot))\goesto H$,
a.s.\ as
$f(n)\goesto\infty$, provided only that $f(n)\leq 
\frac{\log n}{H+\epsilon}$, for some $\epsilon >0$.
It is easy to see that this implies
$\hat{h}_{f(n)}(n)\goesto H$, a.s., for the case
$f(n)=\frac{\log \log n}{\log |A|}$. This proves (a).

To establish part (b), suppose $X\in \cM$ has
order $M$.  The BIC consistency theorem, see 
\cite{Csiszar-Shields}, implies that
	$$
 		 \frac{|A|^{f(n)}(|A|-1)}{2} \log n+n\hat{h}_{f(n)}(n) > 
   \frac{|A|^{M}(|A|-1)}{2}  \log n+n\hat{h}_{M}(n), 
   		$$
eventually a.s.  Using the relation $|A|^{f(n)}=\log n$ and
the   bound $n\hat{h}_{M}(n)\geq nH-c\log\log n$, which
holds eventually a.s.\ by Lemma~\ref{LIL}, we
obtain 
$$
n\hat{h}_{f(n)}(n)\geq nH -c\log \log
n+\frac{|A|^{M}(|A|-1)}{2}  \log n
-\frac{(|A|-1)}{2}\log^2n,
$$
from which (b) follows. \qed

	\begin{remark} 
The empirical entropy estimator $\hat{h}_{f(n)}(n)$
 converges to entropy faster than the
recurrence-based estimator, which is not surprising as the
latter uses so  little about the sample path.
 We suspect there may be
a more direct proof of 
Proposition~\ref{coding-rate-thm}(b) than the one we gave. 
	\label{estimator-remark}       \end{remark} 

	\begin{remark} An important example for
which an $o(1)$ underestimation bound is not known is the
Lempel-Ziv entropy estimator, \cite{Ziv}.  An $O((1/n)\log
n)$ underestimation bound for the class $\cM_0$ of i.i.d.\
processes has been established, see
\cite{Jacquet-Szpankowski}, a result we suspect  can be
extended to the class $\cM$.
\label{LZ}       \end{remark} 

\subsection{The ``flat spot'' problem.}
\hspace{\parindent}\label{flat-spot-sec} For the Markov
order estimation problem, it is tempting to take as order
estimator the first
$k$ for which $\hat{h}_k(n)-\hat{h}_{k+1}(n)<n^{-1/4}$. 
This eventually a.s.\ gets stuck at the first $k$
for which $H_k=H_{k+1}$.  Such flat spots can occur for
$k<M-1$. 
This shows, incidentally, why we needed to take the
maximum over a growing interval of possible orders in
the definition (\ref{testfcn-def}) of our maximal fluctuation test
function.

	\begin{remark} 
The ``no flat spot'' case is  ``generic'' for 
it is easy to see that in the usual parametrization of
the set of $X\in \cM$ of order $M$ as a subset of
$|A|^M(|A|-1)$-dimensional Euclidean space,  the set of
$X$ of order $M$ whose conditional entropy has flat spots
before $M$ has Lebesgue measure 0.  This is a good example
where genericity is not an interesting concept.  
	\label{generic} \end{remark}

\subsection{The BIC, MDL, and related methods.}
\label{related-methods-sec}
\hspace{\parindent}
Two important and related methods, the Bayesian
Information Criterion (BIC) and the Minimum
Description Length (MDL) Principle  are the basis for
many model selection methods, see
\cite{BRY, Csiszar-MDL, Csiszar-Shields} for discussion and
references to these and other methods.  Both the BIC
and the MDL focus on selecting  the correct class
from a nested sequence of parametric model classes,
$\cM_0\subset
\cM_1\subset
\cM_2 \ldots,$ based on a sample path drawn from some
$\mbox{$\PP$}\in\cup\cM_k$.

The BIC, introduced by Schwarz~\cite{Sch}, is based on  Bayesian principles
and  leads to the model estimator 
$$
M^*_{\bic}(x_1^n)\stackrel{\mbox{\scriptsize
def}} =\arg\min_k\biggl(-\log
\mbox{$\PP$}_{\ml(k)}(x_1^n)+
     \frac{\phi(k)}{2}\log n\biggr),
 $$
where $\mbox{$\PP$}_{\ml(k)}(x_1^n)$ is the $k$-th order  maximum
likelihood, i.e., the largest probability given to
$x_1^n$ by distributions in $\cM_k$, and $\phi(k)$ is the number of
free parameters needed to describe members of $\cM_k$.  For the
Markov order estimation problem,
$\cM_k=\{X\in\cM\colon\ M(X)\leq k\}$,
 $-\log
\mbox{$\PP$}_{\ml(k)}(x_1^n)=(n-k)\hat{h}_k(n)$, and
$\phi(k)=|A|^k(|A|-1)$. Schwarz~\cite{Sch} proved consistency if
the model classes are i.i.d.\ exponential families
 and a  bound on the  number of models is
assumed, a result later extended to the Markov case by 
Finesso~\cite{Finesso}.  The first consistency proofs for the
Markov case without an order bound assumption are given in
\cite{Csiszar-Shields}.  The proofs are surprisingly
complicated, though they have been simplified
somewhat in
 \cite{Csiszar-MDL}, which focuses on   MDL consistency.

 The MDL principle, introduced by Rissanen (see \cite{BRY}), is based
 on universal coding ideas. For each $k\leq n$, the sequence $x_1^n$ is encoded
 using a binary code that is ``optimal'' for the class $\cM_k$ and the
 model that has the shortest code length is chosen, that is,
	\begin{equation}
		M^*_{\mdl}(x_1^n) \stackrel{\mbox{\scriptsize
def}} =\arg\min_k \cL_k(x_1^n)
		\label{MDL-estimator}    \end{equation}
where $\cL_k(x_1^n)$ is the length of the  code word
assigned to
$x_1^n$.  Different concepts of ``optimal'' lead to different
estimators. For a discussion of consistency for such estimators
without a prior order bound,
see \cite{Csiszar-Shields} and \cite{Csiszar-MDL}.

%  In the Markov case, one optimality concept  produces
%the BIC, and two other more natural optimality concepts lead to
%slightly different estimators that are equivalent to 
%the BIC if a  bound on the number of models is
%assumed.  These two latter methods have been shown  to be
%inconsistent if no order bound is assumed, but both
%are consistent if the number of models is allowed to
%grow slowly with $n$, see

\noindent{\bf Acknowledgements}.
We are grateful to  I. Kontoyiannis, M. Krishnapur and T. Rud\'as
for helpful comments.

 \end{document}